\documentclass{article}
\usepackage{geometry}
\usepackage{amssymb}
\usepackage{graphics}
\usepackage{epstopdf}

\newtheorem{thm}{Theorem}
\newtheorem{prop}[thm]{Proposition}

\newtheorem{cor}[thm]{Corollary}

\newcommand{\pf}{\noindent {\bf Proof:}\ }
\newcommand{\fp}{\hspace{\fill} $\Box$ \par \bigskip}

\newcommand{\te}{\textstyle}

\title{The Mahler Measure of Parametrizable Polynomials}
\author{Sam Vandervelde}

\begin{document}
\maketitle

\begin{abstract}
Our aim is to explain instances in which the value of the logarithmic Mahler measure $m(P)$ of a polynomial $P\in\mathbb{Z}[x,y]$ can be written in an unexpectedly neat manner.  To this end we examine polynomials defining rational curves, which allows their zero-locus to be parametrized via $x=f(t)$, $y=g(t)$ for $f,g\in\mathbb{C}(t)$.  As an illustration of this phenomenon, we prove the equality
\[ \pi m(y^2+y(x+1)+x^2+x+1) = 2D(i)-{\textstyle\frac{3}{4}}D(\omega), \]
where $\omega=e^{2\pi i/3}$ and $D(z)$ is the Bloch-$\!$Wigner dilogarithm.  As we shall see, formulas of this sort are a consequence of the Galois descent property for Bloch groups.  This principle enables one to explain why the arguments of the dilogarithm function depend only on the points where the rational curve defined by~$P$ intersects the torus $|x|=|y|=1$.  In the process we also present a general method for computing the Mahler measure of any parametrizable polynomial.
\end{abstract}

\section{Introduction}

The purpose of this paper is to explain certain instances in which the logarithmic Mahler measure $m(P)$ of a polynomial $P\in\mathbb{Z}[x,y]$ can be written in an unexpectedly neat manner.  To begin, we offer three examples that illustrate this sort of phenomenon.  In the formulas below $D(z)$ refers to the Bloch-$\!$Wigner dilogarithm, $\omega=e^{2\pi i/3}$, $\xi_5=e^{2\pi i/5}$, $F$ is the field $\mathbb{Q}(\sqrt{-2})$, and $\zeta_F$ is the corresponding zeta function.
\begin{eqnarray} \pi m(y^2+y(x+1)+x^2+x+1) &=& 2D(i)-{\textstyle\frac{3}{4}}D(\omega) \label{exone} \\
\pi m(y^2+y+x^2+x+1) &=& {\textstyle\frac{3}{4}}D(\omega)+{\textstyle\frac{5}{4}}D(\xi_5) -{\textstyle\frac{5}{6}}D(\xi_5^2) \rule{0pt}{1.5em} \label{extwo} \\
5\pi m(y^2+y(x^2+1)+x^4+x^3+x^2+x+1) &=& 9D(\omega)-{\textstyle\frac{4}{3}}D(i) +\frac{16\sqrt{2}}{\pi^2}\:\zeta_F(2) \label{exthree} \end{eqnarray}
All three equalities were conjectured on the basis of numerical evidence; however, we will later provide a proof of the first formula.  As we shall see, these statements are different manifestations of the same underlying principle involving the Galois descent property for Bloch groups.  It is natural at this point to ask for a description of the set of all polynomials to which this principle applies.  We will demonstrate that this set includes (but does not appear to be limited to) certain tempered, parametrizable polynomials.  Along the way we will also present a general method for computing the Mahler measure of any parametrizable polynomial.

We begin with a few preliminary definitions and background.  The logarithmic Mahler measure of a polynomial $P\in\mathbb{C}[x_1,\ldots, x_n]$ is given by
\begin{equation} m(P) = \frac{1}{(2\pi i)^n}\int_{T^n} \log|P(x_1, \ldots, x_n)|
\;\frac{dx_1}{x_1}\cdots\frac{dx_n}{x_n}, \end{equation}
where the $n$-torus $T^n$ consists of all complex $n$-tuples $(x_1, \ldots, x_n)$ with $|x_1|=\cdots=|x_n|=1$.  In particular, this integral converges even though $P(x,y)$ may vanish on~$T^n$.  Despite this somewhat analytical definition, Mahler measure bears surprisingly diverse connections to topics in number theory and geometry.  Thus Boyd~\cite{Bo1} and Rodriguez Villegas~\cite{FRV} have produced polynomials $P\in\mathbb{Z}[x,y]$ for which $m(P)$ may be expressed in terms of the $L$-function of an associated elliptic curve evaluated at $s=2$. Furthermore, Boyd \cite{Bo2} has related Mahler measure to volumes of hyperbolic manifolds, which permits proofs of some identities that he had observed.  Finally, in certain instances Maillot~\cite{Mai}, Vandervelde~\cite{Van}, and Lalin~\cite{Lal} have found ways to express $m(P)$ in terms of the shape, side lengths, and angles within a single cyclic polygon.

In the case of one-variable polynomials it is always possible to obtain an explicit expression for Mahler measure by using Jensen's formula.  One finds that if $P(x)=\lambda(x-\alpha_1)\cdots (x-\alpha_d)$, then \begin{equation} m(P) = \log|\lambda|+\sum_{k=1}^{d} \log^+|\alpha_k|, 
\label{jensen} \end{equation}
where $\log^+ x=0$ for $x<1$ while $\log^+ x=\log x$ when $x\ge1$.  This tool will enable us to perform the first integration for polynomials in two variables.  Clearly $m(PQ)=m(P)+m(Q)$, so we need only analyze irreducible polynomials, which we will assume to be the case in what follows.  Finally, we will often present $P(x,y)$ via its Newton polygon, in which each term $Cx^ly^m$ is represented by labeling the point $(l,m)$ with the coefficient~$C$, as done below.
\[ \begin{array}{rrrrrr}  &2&&&& \\ 
P(x,y)=x^4+4x^3+3x-5-7y+6x^2y+2y^2 \quad\longrightarrow\quad &-7&0&6&& \\
&-5&3&0&4&1  \end{array}. \]
Since multiplication by powers of $x$ or~$y$ does not affect Mahler measure, we will usually neglect to specify the origin, or assume that the Newton polygon `rests' against the axes.  The edge polynomials are obtained by reading off the coefficients, in a counter-clockwise order, along each edge of the polygon.  Thus $P(x,y)$ above has edge polynomials $-5z^4+3z^3+4z+1$, $z^2+6z+2$, and $2z^2-7z-5$.  Following~\cite{FRV}, we say that $P(x,y)$ is \textit{tempered} if all the edge polynomials are cyclotomic; i.e.\ if their roots consist solely of roots of unity.

We next present our objects of study.  We say that a polynomial $P\in\mathbb{C}[x,y]$ is \textit{parametrizable} if it is irreducible and there exist rational functions $f,g\in\mathbb{C}(t)$ which parametrize the curve, i.e.\ provide a bijection (except at singular points) between $\mathbb{P}^1(\mathbb{C})$ and the projective curve defined by~$P$ via the correspondence $t\mapsto (f(t),g(t))$.  In particular, $P(f(t),g(t))$ should be identically zero.  In other words, $P(x,y)$ must define a rational curve.

In the case of parametrizable polynomials, the link between Mahler measure and other areas of mathematics is provided by the Bloch-$\!$Wigner dilogarithm, which we briefly introduce.  It may be defined as
\begin{equation} D(z) = \mathrm{Im}(\mathrm{Li}_2(z)) + \log|z|\arg(1-z), \qquad z\in
\mathbb{C}\smallsetminus[1,\infty). \end{equation}
Unlike the standard dilogarithm function~$\mathrm{Li}_2(z)$, the Bloch-$\!$Wigner dilogarithm (which we shorten to just `the dilogarithm' from here on) may be extended to a function on~$\mathbb{P}^1(\mathbb{C})$, vanishing on the real axis and at infinity.  It is real analytic except at $z=0$, 1, and~$\infty$, where it is merely continuous.  The dilogarithm is invariant under the group of six transformations generated by $z\mapsto 1-\bar{z}$ and $z\mapsto 1-1/z$, and satisfies $D(\bar{z})=-D(z)$ as well.  There is also an analogue to Abel's five-term identity for~$\mathrm{Li}_2(z)$, namely
\begin{equation}
D(a)+D(b)+D\left(\frac{1-b}{a}\right)+ D\left(\frac{a+b-1}{ab}\right)+ D\left(\frac{1-a}{b}\right)=0 
\label{relat} \end{equation}
for any $a,b\in\mathbb{C}^*$.  (See~\cite{Za1} for a more comprehensive introduction to this function.)  The dilogarithm will play a role in evaluating Mahler measure since it appears as the primitive of a certain differential \mbox{1-form}, which will enable us to complete the integration appearing in the definition.

\section{Computational Tools}

Techniques for determining $m(P)$ for particular types of parametrizable polynomials have been developed by various authors, including Smyth~\cite{Smy}, Maillot~\cite{Mai}, and Rodriguez Villegas~\cite{FRV}.  The purpose of this section is to present a systematic way of evaluating $m(P)$ for any such polynomial.  Our approach closely follows the treatment in~\cite{BR2}.

Suppose that $P(x,y)$ has degree~$d$ as a polynomial in~$y$ and that the coefficient of~$y^d$ is~$\lambda(x)$.  When $\lambda(x)\ne 0$ we can write
\[ P(x,y) = \lambda(x)\prod_{k=1}^d (y-\rho_k(x)), \]
where the $\rho_k$ catalog the $d$ roots of $P(x,y)$ for a fixed value of~$x$.  One may then integrate $\log|P(x,y)|$ with respect to $y$ using Jensen's formula:
\[ \frac{1}{2\pi i} \int_{T^1} \!\left( \log|\lambda(x)|+\sum_{k=1}^d \,\log|y-\rho_k(x)| \right) \frac{dy}{y} = \log|\lambda(x)|+\sum_{k=1}^d \,\log^+|\rho_k(x)|. \]
This equality is not valid for the finitely many values of~$x\in T^1$ where $\lambda(x)=0$, but this does not affect our computation of $m(P)$ since we will subsequently be integrating the above expression with respect to~$x$.  According to the definition of Mahler measure we have
\begin{eqnarray} m(P) 
&=& \frac{1}{(2\pi i)^2} \int_{T^2} \log|P(x,y)| \;\frac{dy}{y} \;\frac{dx}{x} \nonumber \\
&=& \frac{1}{2\pi i} \int_{T^1} \!\left(\log|\lambda(x)|+\sum_{k=1}^d \,\log^+|\rho_k(x)|
\right) \frac{dx}{x} \nonumber \\
&=& m(\lambda(x)) + \frac{1}{2\pi} \,\mathrm{Im}\!\int_{T^1} \sum_{k=1}^d \,\log^+|\rho_k(x)| \;\frac{dx}{x}.
\label{xplane}\end{eqnarray}
The last equality stems from the fact that $\log^+|\rho_k(x)|$ is real while $\frac{dx}{x}$ is pure imaginary for $x\in T^1$, as may be seen by taking $x=e^{i\varphi}$.  

Rather than keep track of each of the $d$ roots separately with the functions~$\rho_k(x)$, we adopt the following more unified approach.  Let $\cal C$ be the projective curve defined by $P(x,y)=0$.  (By slight abuse of notation, we will continue to refer to the affine coordinates of points on~$\cal C$.)  Then the set of all points on~$\cal C$ for which $|x|=1$ will consist of one or more closed loops, which we orient compatibly (via the $x$-coordinate) with the usual orientation of~$T^1$.  Clearly the values of the $y$-coordinates are the same as the values of~$\rho_k(x)$.  (Observe that points at infinity on~$\cal C$ occur precisely when $\lambda(x)=0$ for $|x|=1$.)  Now suppose that $P(x,y)$ is parametrizable via $t\mapsto(f(t),g(t))$ for rational functions $f$ and~$g$, and let $\gamma\subset\mathbb{P}^1(\mathbb{C})$ be the oriented pullback of this set.  Then~(\ref{xplane}) may be rewritten as
\[ 2\pi m(P) = 2\pi m(\lambda(x)) + \mathrm{Im}\!\int_{\gamma} \log^+|g(t)| \;\frac{df(t)}{f(t)}. \]
The presence of the invariant measure~$\frac{dx}{x}$ in the definition of Mahler measure is crucial to the equality of these two integrals.

Because $\log^+|g(t)|=0$ when $|g(t)|<1$, we are only interested in values of~$t$ for which $|g(t)|\ge1$.  Therefore let $\gamma_1$, \ldots, $\gamma_n$ be the oriented paths which comprise those portions of~$\gamma$.  We denote the initial and terminal points of $\gamma_j$ by $u_j$ and~$v_j$, respectively.  These endpoints arise from points $(x,y)\in{\cal C}$ for which $|x|=|y|=1$.  Since these points play an important role in what follows,  we will refer to them as the \textit{toric points\/} of~$P(x,y)$.  It is also possible that some of the $\gamma_j$ are loops which lie entirely outside the unit circle, in which case any point on the loop may be chosen to serve as both endpoints.  We can now replace $\log^+$ by $\log$ in the above equation, leading to
\begin{equation} 2\pi m(P) = 2\pi m(\lambda(x)) + \sum_{j=1}^n \,\mathrm{Im}\!\int_{\gamma_j} \log|g| \;\frac{df}{f}. \label{logg} \end{equation}

Following~\cite{BR1} we introduce a closed, real-valued differential 1-form $\eta(f,g)$ on~$\mathbb{P}^1(\mathbb{C})$ for each ordered pair of non-zero rational functions $f,g\in\mathbb{C} (t)$.  It is given by
\[ \eta(f,g) = \log|f|\;d\arg g - \log|g|\;d\arg f, \]
or, equivalently, as
\begin{equation}
\eta(f,g) = \mathrm{Im}\!\left(\log|f|\;\frac{dg}{g}-
\log|g|\;\frac{df}{f}\right). \end{equation}
As may readily be verified, $\eta$~is bimultiplicative and  skew-symmetric in its arguments.  Furthermore, if $|\xi|=1$ then $\eta(\xi f,g)=\eta(f,\xi g)=\eta(f,g)$.  We will use the fact that $\eta(f,1-f)$ is exact with primitive~$D\circ f$, where $D$ is the Bloch-Wigner dilogarithm.

Since $f(t)\in T^1$ for $t\in \gamma_j$, we have $\log|f|=0$ along this path.  Therefore the summand of~(\ref{logg}) can be written as
\begin{equation} \mathrm{Im}\! \int_{\gamma_j} \log|g|  \;\frac{df}{f} = -
\int_{\gamma_j} \eta(f,g). \label{eta} \end{equation}
In order to make further progress we utilize an algebraic technique ascribed to Tate.  It follows from the properties of~$\eta$ that
\[ \eta(t-\alpha,t-\beta) = \eta\left(\frac{t-\alpha}{\beta-\alpha},
1-\frac{t-\alpha}{\beta-\alpha}\right) + \eta(t-\alpha,\alpha-\beta) +
\eta(\beta-\alpha,t-\beta) \]
for~$\alpha\ne\beta$.  Of course, we have $\eta(t-\alpha,t-\alpha)=0$.  Therefore if we write $f(t)=\lambda_1\prod(t-\alpha_r)^{l_r}$ and $g(t)=\lambda_2 \prod(t-\beta_s)^{m_s}$ for non-zero integers $l_r$ and~$m_s$ and $\lambda_1, \lambda_2\in\mathbb{C}$, then
\begin{eqnarray*} \eta(f,g) &=& \sum_{r,s} l_rm_s\eta(t-\alpha_r,
t-\beta_s)+\sum_r l_r\eta(t-\alpha_r,\lambda_2) + \sum_s m_s
\eta(\lambda_1,t-\beta_s) \\ 
&=& {\sum_{r,s}}' l_rm_s\eta\left(\frac{t-\alpha_r}{\beta_s-\alpha_r},
1-\frac{t-\alpha_r}{\beta_s-\alpha_r}\right) 
+ {\sum_{r,s}}' l_rm_s\eta(t-\alpha_r,\alpha_r-\beta_s) + {}\\
&& {\sum_{r,s}}' l_rm_s\eta(\beta_s-\alpha_r,t-\beta_s) +\sum_r
l_r\eta(t-\alpha_r,\lambda_2) + \sum_s m_s \eta(\lambda_1,t-\beta_s).
\end{eqnarray*}
The notation ${\sum}'$ indicates the sum over all pairs $r$ and~$s$ where $\alpha_r\ne\beta_s$.  We also used $\eta(\lambda_1,\lambda_2)=0$ in the first line.  Collecting terms involving $(t-\alpha_r)$ yields $\sum_r l_r\eta(t-\alpha_r, \lambda_2\prod'_s (\alpha_r-\beta_s)^{m_s})$.  Let us define $\tilde{g}$ by $\tilde{g}(t)=g(t)$ unless $t=\beta_s$ for some~$s$, in which case the corresponding factor in
$g$ is omitted before evaluating at~$t$.  Hence $\tilde{g}(\alpha_r) \ne 0, \infty$.  Combining terms involving $(t-\beta_s)$ and defining $\tilde{f}$ similarly, this algebra can be summarized as follows.
\begin{prop} \label{tatethm} With notation as above, 
\begin{eqnarray}
 \eta(f,g) &=& {\sum_{r,s}}' l_rm_s\eta\left(\frac{t-\alpha_r}
{\beta_s-\alpha_r}, 1-\frac{t-\alpha_r}{\beta_s-\alpha_r}\right) + {}
\nonumber \\
&& \sum_r l_r\eta(t-\alpha_r,\tilde{g}(\alpha_r)) 
+ \sum_s m_s\eta(\tilde{f}(\beta_s),t-\beta_s). 
\label{tate} \end{eqnarray} \end{prop}

Recall that the endpoints of $\gamma_j$ are $u_j$ and~$v_j$.  Using the fact that $\eta(f,1-f)$ is exact with primitive $D\circ f$ and invoking Stoke's theorem, the first term of~(\ref{tate}) integrates to 
\[ {\sum_{r,s}}' l_rm_s \left( D\left(\frac{v_j-\alpha_r}
{\beta_s-\alpha_r} \right) - D\left(\frac{u_j-\alpha_r}
{\beta_s-\alpha_r} \right) \right). \]
According to the definition of~$\eta$, the second sum in~(\ref{tate}) will become
\[ -\sum_r l_r \log|\tilde{g}(\alpha_r)| \;\mathrm{Im}\! \int_{\gamma_j} \frac{dt}{t-\alpha_r}. \]
The integral appearing here calculates the ``winding angle'' of the path~$\gamma_j$ around the point~$\alpha_r$, so we employ the notation
\begin{equation} \mathrm{wind}(\gamma_j,\alpha_r) = \mathrm{Im}\!\int_{\gamma_j} \frac{dt}{t-\alpha_r}. \label{wind} \end{equation}
Observe that since $\gamma_j$ is differentiable, this expression has a limiting value even if the path passes through either $\alpha_r$ or the point at infinity.  In fact, $\alpha_r$ cannot lie on~$\gamma_j$, because $\alpha_r$ is either a zero or pole of~$f$, while the image of~$\gamma_j$ under~$f$ is a subset of the unit circle.  In the same manner, winding angles will appear upon integrating the third sum in~(\ref{tate}).  It is possible for some point~$\beta_s$ to lie on~$\gamma_j$; however, in this case $f(\beta_s)\in T^1$, since $\beta_s\in\gamma_j$.  In particular, $\beta_s$~will not be a zero or pole of~$f(t)$, so $\tilde{f}(\beta_s)=f(\beta_s)$ and $\log|\tilde{f}(\beta_s)|=0$.  Therefore the corresponding term $\log|\tilde{f}(\beta_s)|\cdot\mathrm{wind}(\gamma_j,\beta_s)$ arising from the third sum in~(\ref{tate}) will vanish.

By combining the preceding discussion with~(\ref{eta}) we are able to evaluate~$m(P)$ in terms of standard functions.  We summarize our findings, utilizing the above notation.
\begin{thm} Suppose that $P\in\mathbb{C}[x,y]$ can be parametrized by $x=f(t)=\lambda_1\prod(t-\alpha_r)^{l_r}$ and $y=g(t)=\lambda_2 \prod(t-\beta_s)^{m_s}$.  Let~$S$ consist of those points in the zero-locus of~$P$ for which $|x|=1$ and $|y|\ge 1$, and let $\gamma_1$, \ldots, $\gamma_n$ be the paths which map to~$S$ under $t\mapsto(f(t),g(t))$, oriented via $f(\gamma_j)\subset T^1$.  If $u_j$ and~$v_j$ denote the initial and terminal points of~$\gamma_j$, and $P(x,y)$ has leading coefficient $\lambda(x)$ as a polynomial in~$y$, then
\begin{eqnarray} 2\pi m(P) &=& 2\pi m(\lambda(x)) + \sum_{j=1}^n \left[{\sum_{r,s}}'
l_rm_s \left( D\left(\frac{u_j-\alpha_r}{\beta_s-\alpha_r} \right) - D\left(\frac{v_j-\alpha_r}
{\beta_s-\alpha_r} \right) \right) + {} \right. \nonumber \\
&& \left. \sum_r l_r \log|\tilde{g}(\alpha_r)| \cdot \mathrm{wind}
(\gamma_j,\alpha_r) - \sum_s m_s \log|\tilde{f}(\beta_s)| \cdot \mathrm{wind}
(\gamma_j,\beta_s) \right].
\label{tool}\end{eqnarray}
\label{eval}\end{thm}

\section{Applying the Evaluation Theorem}

\begin{figure} \centerline{\includegraphics{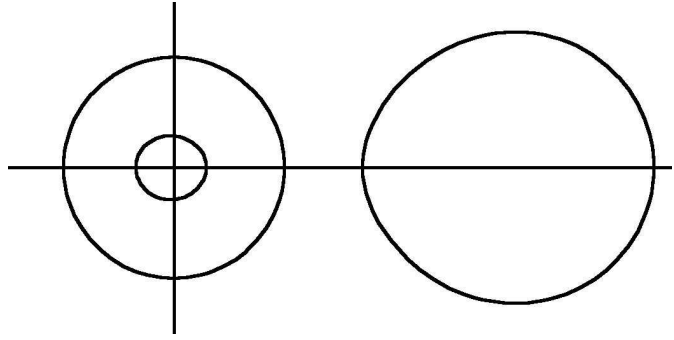}}
\caption{The values of $y$ satisfying $P(e^{i\varphi},y)=0$ for $P(x,y)=-2y^2+2xy+6y+2x+1$.}
\label{graph1} \end{figure}

We next present three examples to highlight some of the different considerations and algebraic methods that arise when computing Mahler measure.  To begin, we take 
\[ P(x,y)=-2y^2+2xy+6y+2x+1, \] which has the obvious parametrization 
\[ y=g(t)=t, \qquad x=f(t)=\frac{2t^2-6t-1}{2t+2} = \frac{(t-\frac{1}{2}(3+\sqrt{11})) (t-\frac{1}{2}(3-\sqrt{11}))} {t-(-1)}. \]
Setting $x=e^{i\varphi}$, solving for $y$ in $P(x,y)=0$, and then plotting the result yields the graph shown in Fig.~\ref{graph1}.  (The unit circle is also included for reference.)  
The loop outside the unit circle is the set~$S$ alluded to in Theorem~\ref{eval}.  Since $y=t$, it is also the single path $\gamma_1$ that maps to~$S$.  This loop is oriented in a counterclockwise manner, intersecting the real axis at $t=\frac{1}{2}(2+\sqrt{2})$ and $\frac{1}{2}(4+\sqrt{22})$.  Because $u_1=v_1$ for a loop, the actual value we choose does not matter; the dilogarithm terms in~(\ref{tool}) will cancel.  The only zero or pole of either $f(t)$ or $g(t)$ that lies within the loop is $\alpha_1=\frac{1}{2}(3+\sqrt{11})$, hence $\mathrm{wind}(\gamma_1,\alpha_1)=2\pi$ while all other winding angles equal 0.  Since $\tilde{g}(\alpha_1)=\frac{1}{2}(3+\sqrt{11})$ and $\lambda(x)=-2$, we conclude that
\begin{equation} m(P) = \log(2) + \log({\textstyle\frac{1}{2}}(3+\sqrt{11})) = \log(3+\sqrt{11}),
\end{equation}
which can be confirmed numerically.  This example demonstrates the general principle that if $P(x,y)$ has no toric points, then $m(P)$ may be expressed in terms of logarithms of algebraic integers.

As we have just seen, finding a parametrization is trivial when $P(x,y)$ is linear in either $x$ or~$y$.  There are also standard techniques for handling degree two polynomials.  These have Newton polygons in the shape of a small triangle:
\[ \begin{array}{cccc}  &c&& \\ 
P(x,y)=ax^2+bxy+cy^2+dx+ey+f\quad\longrightarrow\quad &e&b& \\
&f&d&a  \end{array}. \]
Recall that the edge polynomials are $fz^2+dz+a$, $az^2+bz+c$, and $cz^2+ez+f$.  We present one parametrization which is of particular relevance.
\begin{prop} Let $P\in\mathbb{Z}[x,y]$ be an irreducible degree two polynomial, and let $E$ be the field obtained by adjoining to~$\mathbb{Q}$ the roots of all the edge polynomials of~$P$.  Then there exist rational functions $f(t)$ and~$g(t)$ which parametrize~$P$, all of whose zeros, poles, and coefficients lie in~$E$. \label{param} \end{prop}
\pf  The existence of such a parametrization is something of an algebraic \textit{fait accompli.}  Write $\Delta=b^2-4ac$, $\Delta_x=e^2-4cf$, and $\Delta_y=d^2-4af$, the discriminants of the edge polynomials, and choose $\kappa\in E$ satisfying $c\kappa^2+b\kappa+a=0$.  Then if $\Delta\ne0$ we claim that
\begin{eqnarray} f(t)&=&\frac{1}{\Delta}\cdot\frac{t^2+(2cd-be)t+c(ae^2+
fb^2+cd^2-bde-4acf)}{t}, \rule{2cm}{0pt} \nonumber \\
g(t)&=&\frac{1}{\Delta}\cdot\frac{(\kappa t)^2+(2ae-bd)
(\kappa t)+a(ae^2+fb^2+ cd^2-bde-4acf)}{\kappa t} \rule{0pt}{1cm}
\label{canon}\end{eqnarray}
is a parametrization for~$P(x,y)$.  The identity $P(f(t),g(t))\equiv 0$ is easily verified using a computer algebra system.  We also observe that the discriminants of the numerators of~$f(t)$ and~$g(t)$ are $\Delta\cdot \Delta_x$ and $\kappa^2\Delta\cdot \Delta_y$, so their roots lie in~$E$, as desired.

If $\Delta=0$ our polynomial must be of the form $P(x,y)=c'(a'x+b'y)^2+dx+ey+f$.  We then define $\Delta'=c'(a'e-b'd)$ and take
\begin{equation} f(t)=\frac{1}{\Delta'}\,(b't^2+et+b'c'f), \qquad
g(t)=-\frac{1}{\Delta'}\,(a't^2+dt+a'c'f). \end{equation}
In this case the discriminants are simply $\Delta_x$ and~$\Delta_y$, respectively.  Finally, we note that if $\Delta'=0$, then $P(x,y)$ is a function of $(a'x+b'y)$, so that $P(x,y)$ factors over~$\mathbb{C}$, contrary to hypothesis. \fp

\begin{figure} \centerline{\includegraphics{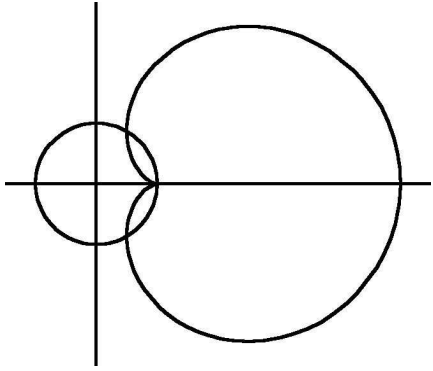}}
\caption{The values of $y$ satisfying $P(e^{i\varphi},y)=0$ for $P(x,y)=x^2-2xy+y^2-4y+4$.}
\label{graph2} \end{figure}

In this manner we obtain a parametrization for $P(x,y)=x^2-2xy+y^2-4y+4$, our second example.  Since $\Delta=0$ here, we use the alternate formulas with $\Delta'=-4$ to find $f(t)=\frac{1}{4}(t+2)^2$ and $g(t)=\frac{1}{4}(t^2+4)$.  Substituting $2t$ for~$t$ yields the slightly nicer functions $f(t)=(t+1)^2$ and $g(t)=t^2+1$.  The lead coefficient of $P(x,y)$ is $\lambda(x)=1$, so $m(\lambda(x))$ contributes nothing in~(\ref{tool}).  The set of points in the zero-locus of~$P$ for which $|x|=1$ is a single loop; its projection onto the $y$-coordinate is shown in Fig.~\ref{graph2}.
It is not hard to determine that the path $\gamma_1$ mapping to the portion of this loop outside the unit circle is described by $t=e^{i\varphi}-1$ for $\frac{\pi}{3} \le\varphi\le \frac{5\pi}{3}$.  This path has initial point $u_1=\omega$ and terminal point $v_1=\overline{\omega}$, where $\omega=e^{2\pi i/3}$.  The zeros and poles of~$f(t)$ and~$g(t)$ are $\alpha_1=-1$ with multiplicity $l_1=2$, and $\beta_1=i$, $\beta_2=-i$ with multiplicities $m_1=m_2=1$.  The winding angles of $\gamma_1$ about these
points are easily found to be
\[ \mathrm{wind}(\gamma_1,-1)=\frac{4\pi}{3}, \quad
\mathrm{wind}(\gamma_1,i)=\frac{\pi}{3}, \quad\mbox{and}\quad
\mathrm{wind}(\gamma_1,-i)=\frac{\pi}{3}. \]
Therefore our evaluation theorem predicts that
\begin{eqnarray*} 2\pi m(P) &=&
2\left(D\left(\frac{\omega+1}{i+1}\right)
-D\left(\frac{\overline{\omega}+1}{i+1}\right)
+D\left(\frac{\omega+1}{-i+1}\right) 
-D\left(\frac{\overline{\omega}+1}{-i+1}\right)\right) \\
&& + 2\log 2\cdot\frac{4\pi}{3} - \log 2\cdot\frac{\pi}{3} -
\log 2\cdot\frac{\pi}{3}.
\end{eqnarray*}
Although it is not immediately obvious, this expression may be reduced to
\begin{equation} 
\pi m(P) = \frac{15}{4}\: D(\omega) + \pi \log 2. \end{equation}
Notice that the coordinates of the toric points are present in the simplified expression, but the roots of the edge polynomials have disappeared.

We pause to mention that a toric point $(x,y)$ will also satisfy the polynomial $x^ly^m P(\frac{1}{x},\frac{1}{y})$, where $l$ and~$m$ are the degrees of $x$ and~$y$ in~$P$.  This follows from conjugating $P(x,y)=0$ and then using $\overline{x}=1/x$ and $\overline{y}=1/y$.  Taking the resultant with respect to~$y$ yields a polynomial in~$x$ whose roots include all the desired $x$-coordinates.  Applying this method to the above example, one finds that the toric points are $(\omega,-\omega)$, $(\bar{\omega},-\bar{\omega})$ and $(1,1)$, which agrees with Fig.~\ref{graph2}.  Note that this technique breaks down in the special case of reciprocal polynomials, where $x^ly^mP(\overline{x},\overline{y})$ is identical to $P(x,y)$.  This occurs when the Newton polygon is centrally symmetric, as in our final example.

To conclude this section we compute $m(P)$ for the polynomial
\[ P(x,y) = x + y - 4xy + x^2y + xy^2. \]
In~\cite{Bo1}, Boyd employs a clever technique which shows that $\pi m(P) = 4D(i)$. We will verify this equality using the machinery developed above.  Polynomials of this form typically define elliptic curves, but this particular example has a singularity at $(1,1)$ and can be parametrized by
\[ f(t) = \frac{(t-1)(t-i)}{(t+1)(t+i)}, \quad g(t) = \frac{(t-1)(t+i)}{(t+1)(t-i)}.  \]
The $y$-values appearing in the set~$S$ range from 1 out to $3+2\sqrt{2}$ and back along the real axis, and the path~$\gamma_1$ mapping to~$S$ is the ray consisting of the positive real multiples of 
$(-1+i)$, with initial point $u_1=\infty$ and terminal point $v_1=0$.  Simple geometry dictates that the winding angles are
\[ \mathrm{wind}(\gamma_1,1) = \frac{\pi}{4}, \quad \mathrm{wind}(\gamma_1,i) = \frac{3\pi}{4}, \quad \mathrm{wind}(\gamma_1,-1) = -\frac{3\pi}{4}, \quad \mathrm{wind}(\gamma_1,-i) = -\frac{\pi}{4}. \]
The values of $|\tilde{f}(\beta_s)|$ and $|\tilde{g}(\alpha_r)|$ all reduce to either 2 or~$\frac{1}{2}$, and one finds that the angle terms in~(\ref{tool}) all cancel in pairs.  Using the fact that $D(\infty)=0$, the dilogarithm terms become
\[ 2\pi m(P) = 2D\!\left(\frac{1}{1-i}\right) + 2D\!\left(\frac{i}{i+1}\right) - 2D\!\left(\frac{1}{1+i}\right) - 2D\!\left(\frac{i}{i-1}\right). \]
Dividing by~2, rationalizing denominators, and utilizing $D(z)=-D(\bar{z})=D(1-\frac{1}{z})$ yields 
\begin{equation} \pi m(P)=4D\!\left(\frac{1+i}{2}\right)=4D(i). \end{equation}

\section{The Bloch Group and Galois Descent}

Our primary tool for addressing the dilogarithm terms appearing in~(\ref{tool}) is the Galois descent property for Bloch groups.  We begin by describing (one version of) the Bloch group of a field~$K$.  For any element $\theta=\sum c_k[\alpha_k]\in\mathbb{Z}[K]$, we define $\partial(\theta)\in\wedge^2 K^*$ by
\begin{equation} \partial(\theta) 
= \sum_{k=1}^n c_k\; \alpha_k\wedge(1-\alpha_k), \end{equation}
where we omit any term of the sum corresponding to $\alpha_k=0$ or~1.  The kernel of $\partial\!:\!\mathbb{Z}[K]\rightarrow\wedge^2 K^*$ is a subgroup $B'(K)$ of $\mathbb{Z}[K]$ containing all elements of the form
\[ [\alpha]+[\beta]+\left[\frac{1-\beta}{\alpha}\right]+\left[\frac
{\alpha+\beta-1}{\alpha\beta}\right]+\left[\frac{1-\alpha}{\beta}\right] \]
for $\alpha, \beta\in K^*$.  Let $B''(K)$ be the subgroup of~$B'(K)$ generated by these elements.  Then the Bloch group $B(K)$ is defined as the quotient $B'(K)/B''(K)$.  Now assume that $K\subset\mathbb{C}$ via a particular complex embedding.  By linearity we may extend the dilogarithm to $\mathbb{Z}[K]$ via $D(\theta)=\sum c_k D(\alpha_k)$.  In light of the five-term identity~(\ref{relat}), the dilogarithm vanishes on~$B''(K)$, so we obtain a well-defined homomorphism $D\!:\!B(K) \rightarrow \mathbb{R}$.

The Galois descent property for Bloch groups asserts that, up to torsion, elements of~$B(K)$ may be treated as elements of a smaller Bloch group~$B(F)$, where $F$ is a subfield of~$K$, if they are fixed by all the automorphisms in $\mathrm{Gal}(K/F)$.  Our statement is based upon the one appearing in~\cite{BR2}.
\begin{thm} Let $F\subset K\subset \mathbb{C}$ be number fields with $K$ a Galois extension of~$F$ and $G=\mathrm{Gal}(K/F)$.  For $\sigma \in G$ and $\theta\in B(K)$, we have $\theta^\sigma
=\sum c_k[(\alpha_k)^{\sigma}] \in B(K)$ also.  We write $B(K)^G$ to denote the set of all $\theta\in B(K)$ which are fixed by each $\sigma\in G$.  Then
\begin{equation} B(F)\otimes_{\mathbb{Z}}\mathbb{Q} =
B(K)^G\otimes_{\mathbb{Z}}\mathbb{Q}. \end{equation} \end{thm}
Since the dilogarithm vanishes on torsion elements of~$B(K)$, the Galois descent property implies that for $\theta\in B(K)^G$, we have $D(\theta)=cD(\theta')$ for some $\theta'\in B(F)$ and $c\in\mathbb{Q}$.

Our next task is to pinpoint the fields in which $u_j$ and~$v_j$ reside, the endpoints of the path $\gamma_j$ from Theorem~\ref{eval}.  Given $P\in\mathbb{Z}[x,y]$, let $E$ be the number field obtained by adjoining to~$\mathbb{Q}$ the roots of all the edge polynomials of~$P(x,y)$.  We first remark that when $P$ is parametrizable, there always exist suitable functions $f(t)$ and~$g(t)$, all of whose zeros, poles, and coefficients lie in~$E$.  If $P$ is linear in one of its variables then the claim is obvious, and if $P$ has degree two it follows from Proposition~\ref{param}.  (In fact, any non-singular parametrizable polynomial may be reduced to one of these two cases by operations which preserve Mahler measure.)  However, a standard argument from algebraic geometry demonstrates that such a parametrization exists even when $P$ is singular.  Strictly speaking, the arguments in our main result are applicable as long as $f,g\in E(t)$ for some Galois extension $E/\mathbb{Q}$, which clearly can be arranged.

We now associate a number field $F$ to each algebraic point $(\mu,\nu)$ in the zero-locus of~$P$.  Consider the related polynomial $P_o(x,y)=P(x+\mu,y+\nu)$.  Since $P_o(0,0)=0$, this polynomial has no constant term, so in its Newton polygon one of the edges will ``face the origin.''  As an illustration, consider the Newton polygon below, where `$\circ$' indicates the origin.
\renewcommand{\arraystretch}{.8}
\[ \begin{array}{rrrrrr}  3&1&6&&& \\ 2&2&1&5&& \\ &&2&4&3&3 \\ \circ &&&&1&8
\end{array} \] \renewcommand{\arraystretch}{1}%
In this case the edge with polynomial $2z^2+2z+1$ faces the origin.  Now adjoin $\mu$, $\nu$, and the roots of this edge polynomial to $\mathbb{Q}$ to obtain the desired field~$F$.  Note that if $(\mu,\nu)$ is a non-singular point of~$P$, then the $x$ and~$y$ terms of $P_o$ cannot both vanish, so there will be only a single linear edge polynomial.  Hence $F=\mathbb{Q}(\mu,\nu)$ in this case.

In the subsequent arguments we will make the simplifying assumption that there is only a single edge facing the origin and that its polynomial is irreducible over $\mathbb{Q}(\mu,\nu)$.  (In other words, $P$~is sufficiently well-behaved at its singular points.)  In theory, it would be possible to associate several fields to an exceptional point, one for each irreducible factor of each edge polynomial facing the origin.  (Even more care is required when the edge polynomials contain repeated roots.)  The ensuing theory would still hold, but the bookkeeping would be needlessly complicated, since the above assumption is sufficient in practice.

Recall that a particular $u_j$ (or~$v_j$) is the solution for~$t$ to a pair of equations of the form $f(t)=\mu$, $g(t)=\nu$, where $(\mu,\nu)$ is a toric (and hence algebraic) point of $P$ and $f,g\in E(t)$.  The following result shows that $u_j$ belongs to a field that is independent of the particular parametrization chosen.
\begin{prop} Suppose that $P\in\mathbb{Z}[x,y]$ may be parametrized by $f,g\in E(t)$ for a number field~$E$.  Let $(\mu,\nu)$ be an algebraic point in the zero-locus of~$P$ with an associated field $F$ as described above.  If $\rho$ maps to $(\mu,\nu)$ via $t\mapsto(f(t),g(t))$, then $\rho\in K=EF$. \label{edgepoly} \end{prop}

\pf It will be useful to define the field $L=E(\mu,\nu)$, so that $L\subseteq K$.  By construction, $P_o(x,y)=P(x+\mu,y+\nu)$ provides a polynomial relation between $f_o(t)=f(t)-\mu$ and $g_o(t)=g(t)-\nu$.  Note also that $f_o, g_o\in L(t)$ with $f_o(\rho)=g_o(\rho)=0$, so $\rho$ is a common zero of $f_o$ and~$g_o$, say with multiplicities $l$ and~$m$.  Let $h_o(t)\in L[t]$ be the monic, irreducible polynomial  having $\rho$ as a root, so that we can write $f_o=h_o^l\tilde{f}_o$ and $g_o=h_o^m\tilde{g}_o$ for some $\tilde{f}_o, \tilde{g}_o\in L(t)$.  (We employ the tilde in the same spirit as before, but with a slightly different meaning here.)  Since $P_o$ has no constant term, every monomial term $C f_o(t)^rg_o(t)^s$ in $P_o(f_o,g_o)$ will involve a power of $h_o(t)$, with exponent $lr+ms>0$.  Let~$M$ be the minimal such power occurring, and consider 
\begin{equation} \left. \frac{P_o(f_o(t),g_o(t))}{h_o(t)^M} \right|_{t=\rho}. \label{diffM} 
\end{equation}
Only the term(s) involving the minimal power $h_o(t)^M$ will remain after substituting $t=\rho$.  Each such term will contribute $C\tilde{f}_o(\rho)^r \tilde{g}_o(\rho)^s$.  On the other hand, the sum of these terms equals zero since $P_o(f_o,g_o)\equiv 0$.  We deduce that there are at least two such terms, since each is non-zero.  Hence the monomials of the form $C x^ry^s$ with $lr+ms=M$ constitute the edge polynomial $p_o(z)$ facing the origin in the Newton polygon of~$P_o$.  Dividing through by the highest power of $\tilde{f}_o(\rho)$ in~(\ref{diffM}) reveals that $\tilde{g}_o^l(\rho)/\tilde{f}_o^m(\rho)$ is a root of~$p_o(z)$.  To ease notation, we write $q(t)=\tilde{g}_o^l(t)/\tilde{f}_o^m(t)\in L(t)$.

Let $\rho=\rho_1$, $\rho_2$, \ldots, $\rho_d$ be the distinct roots of~$h_o(t)$.  Employing the above reasoning, we see that each $q(\rho_k)$ is a root of~$p_o(z)$.  We also claim that these numbers are distinct.  For suppose that $q(\rho_1) = q(\rho_2)$.  Replacing each factor of $(t-\rho_2)$ in $f_o$ and
$g_o$ by a small perturbation $(t-(\rho_2+\epsilon))$ yields rational functions $f_{\epsilon}, g_{\epsilon}\in\mathbb{C}(t)$ which are related by a polynomial $P_{\epsilon}\in\mathbb{C}[x,y]$ having the same terms as $P_o$ but slightly different coefficients.  Denote by $p _{\epsilon}(z)$ the edge polynomial facing the origin.  Letting $\epsilon\rightarrow 0$ through values for which $q(\rho_1)$ and~$q(\rho_2+\epsilon)$ are unequal, we see that $p_{\epsilon}(z)$ contains two distinct roots which become a double root of $p_o(z)$ in the limit.  But this contradicts the fact that $p_o(z)$ is irreducible and hence has only simple roots.

The previous discussion implies that the degree of $p_o(z)$ is at least as large as the degree of~$h_o(t)$.  Since $\rho$ is a root of $h_o(t)$ and $q(\rho)=\tau$ is a root of~$p_o(z)$, this means that
 $[L(\tau):L]\ge[L(\rho):L]$.   However, $\tau=q(\rho)\in L(\rho)$, which forces $L(\tau) = L(\rho)$.  Therefore $\rho\in L(\tau)\subset K$, as desired. \fp

In order to apply the Galois descent property for Bloch groups we must relate the dilogarithm terms appearing in~(\ref{tool}) to an element of some Bloch group.  This can be done by imposing the following two conditions on our polynomial.  First, we require that $P(x,y)$ be tempered, meaning that its edge polynomials are all cyclotomic.  Secondly, we need the coordinates of each toric point $(\mu,\nu)$ to be commensurate, in the sense that $\mu^l=\nu^m$ for some pair of integers $(l,m)\ne(0,0)$ to ensure that $lm(\mu\wedge\nu)=0$.  (These two conditions appear to be independent; the second is typically satisfied because either $\mu=\nu$, $\mu=\nu^{-1}$, or $\mu$~and~$\nu$ are roots of unity.)  Recall that we are also assuming that $P(x,y)$ is reasonably well-behaved at any singular toric points, as described above.  A polynomial satisfying these restrictions will be denoted as an \textit{admissible\/} polynomial.

\begin{prop} Let $P\in\mathbb{Z}[x,y]$ be admissible, parametrized by $x=f(t)=\lambda_1\prod(t-\alpha_r)^{l_r}$ and $y=g(t)=\lambda_2 \prod(t-\beta_s)^{m_s}$, with $\lambda_i, \alpha_r, \beta_s\in E$, the field obtained by adjoining the roots of all the edge polynomials to~$\mathbb{Q}$.  Suppose that $(\mu,\nu)=(f(u),g(u))$ is a toric point of~$P$ with associated field~$F$.  Writing $K=EF$ and
\[ \theta = {\sum_{r,s}}' l_rm_s\left[\frac{u-\alpha_r}{\beta_s-\alpha_r}\right] \in \mathbb{Z}[K], \]
then $N\theta$ is an element of the Bloch group $B(K)$ for some positive integer~$N$.  \label{bloch}\end{prop}

\pf To begin, $u\in K$ by Proposition~\ref{edgepoly}.  Also note that $u\ne \alpha_r, \beta_s$ since $|f(u)|=|g(u)|=1$.  Hence the sum defining $\theta$ will not involve a term $[0]$ or~$[1]$.  Since the wedge product $f\wedge g$ and the differential 1-form $\eta(f,g)$ share similar properties, the same algebra leading up to Proposition~\ref{tatethm} implies that
\begin{eqnarray}
 f(u)\wedge g(u) &=& {\sum_{r,s}}' l_rm_s\left(\frac{u-\alpha_r}
{\beta_s-\alpha_r}\right)\wedge \left(1-\frac{u-\alpha_r}{\beta_s-\alpha_r}\right) + {}
\nonumber \\
&& \sum_r l_r(u-\alpha_r)\wedge \tilde{g}(\alpha_r)
- \sum_s m_s(u-\beta_s)\wedge \tilde{f}(\beta_s). \label{wedge}\end{eqnarray}
By assumption, $N_1 f(u)\wedge g(u) = N_1 \mu\wedge\nu = 0$ for some positive integer~$N_1$.  Furthermore, since $P$ is tempered, there is a positive integer $N_2$ such that $\xi^{N_2}=1$ for any root $\xi$ of an edge polynomial of~$P$.  We claim that
\[ N_2 \sum_r l_r(u-\alpha_r)\wedge \tilde{g}(\alpha_r)
- N_2 \sum_s m_s(u-\beta_s)\wedge \tilde{f}(\beta_s) = 0. \]
Suppose that a given $\alpha_r$ is a root of $f(t)$ but is not a zero or pole of~$g(t)$, so that $\tilde{g}(\alpha_r)=g(\alpha_r)$.  Letting $t=\alpha_r$ in the identity $P(f,g)\equiv 0$ yields $P(0,g(\alpha_r))=0$, which implies that $g(\alpha_r)$ is the root of an edge polynomial.  (The leftmost vertical one, to be precise.)  Hence
\[ N_2\ l_r (u-\alpha_r)\wedge \tilde{g}(\alpha_r) = l_r (u-\alpha_r)\wedge g(\alpha_r)^{N_2} = l_r (u-\alpha_r)\wedge 1 = 0. \]
If $\alpha_r$ happens to be a pole of $f(t)$ instead then we examine $Q(x,y)=x^l P(x^{-1},y)$, where $l$ is the degree of~$x$ in~$P(x,y)$.  But $Q(x,y)$ may be parametrized by $x=f(t)^{-1}$, $y=g(t)$, so the same reasoning implies that $g(\alpha_r)$ is a root of the leftmost vertical edge polynomial of~$Q$, which is the rightmost vertical edge polynomial of~$P$.  Similar considerations apply to $\tilde{f}(\beta_s)$ when $\beta_s$ is neither a root nor a pole of~$f(t)$.

It remains to handle the instance where $\alpha_r=\beta_s$ for some $r$ and~$s$.  Suppose that this common value is a zero of both $f(t)$ and~$g(t)$.  Then the same argument used in Proposition~\ref{edgepoly} reveals that $\tilde{g}(\alpha_r)^{l_r}/\tilde{f}(\beta_s)^{m_s}$ is a root of an edge polynomial of~$P$ facing the origin.  Therefore
\[ N_2\ l_r(u-\alpha_r)\wedge \tilde{g}(\alpha_r)
- N_2\ m_s(u-\beta_s)\wedge \tilde{f}(\beta_s) = N_2\ (u-\alpha_r)\wedge 
\frac{\tilde{g}(\alpha_r)^{l_r}}{\tilde{f}(\beta_s)^{m_s}} = 0 \]
as before.  By considering a related polynomial $Q(x,y)$ as above we deduce that these terms cancel regardless of whether $\alpha_r=\beta_s$ is a zero or pole of $f(t)$ and~$g(t)$.  In summary, the only remaining terms of~(\ref{wedge}) are those corresponding to~$\partial(\theta)$.  We conclude that $N\partial(\theta)=0$ for $N=N_1N_2$, so~$N\theta\in B(K)$. \fp

Since $E/\mathbb{Q}$ is Galois, it follows that $K/F$ is also a Galois extension, since $K=EF$.  Thus we are now in a position to apply the Galois descent property for Bloch groups to the element~$N\theta$.  Unfortunately, it is usually the case that $N\theta$ is not fixed by every automorphism $\sigma\in \mathrm{Gal}(K/F)$.  This difficulty is circumvented by considering an entire collection of parametrizations together, allowing us to prove our main result.

\begin{thm} Let $P\in\mathbb{Z}[x,y]$ be an admissible polynomial.  Denote the toric points of~$P$ by $(\mu_j,\nu_j)$ with associated fields~$F_j$.  Then
\begin{equation} \pi m(P) = \sum c_j D(\vartheta_j), \qquad \vartheta_j\in B(F_j),\ \ c_j\in \mathbb{Q^*}.
\end{equation} \end{thm}

Note that this formula depends only on quantities intrinsic to~$P$, independent of the parametrization.  Before turning to its proof we mention a couple of corollaries which explain the conjectured equalities given at the outset of this article.  In general, the dimension of $B(F_j)\otimes\mathbb{Q}$ as a $\mathbb{Q}$-vector space is equal to the number of pairs of complex embeddings of~$F_j$.  For instance, if $F_j=\mathbb{Q}(e^{2\pi i/m})$ is a cyclotomic field, then a basis for $B(F_j)\otimes\mathbb{Q}$ is given by the $\varphi(m)/2$ primitive $m^{\mathrm{th}}$ roots of unity lying in the upper half plane.  We therefore obtain

\begin{cor} Let $P\in\mathbb{Z}[x,y]$ be an admissible polynomial.  If the fields $F_j$ associated with the toric points are cyclotomic, then
\[ \pi m(P) = \sum c_j D(\xi_j), \]
where the $c_j$ are rational numbers and the sum ranges over all primitive roots of unity $\xi_j$ in the upper half plane contained in the fields~$F_j$. \label{corone} \end{cor}

On the other hand, Mahler measure may also be related to special values of zeta-functions as a consequence of the following theorem of Borel, found in~\cite{Za2}.

\begin{thm} \label{borelthm}
Let $F$ be a number field with one pair of complex embeddings and $r$ real embeddings, and take $F\subset \mathbb{C}$ via one of the complex embeddings.  If $\theta\in B(F)$ is a non-torsion element, then
\begin{equation} D(\theta) = c\,\frac{|\mathrm{disc}(F)|^{3/2}}
{\pi^{2r+2}} \zeta_F(2), \qquad c\in\mathbb{Q}^*.
\label{borel} \end{equation} \end{thm}

\begin{cor} Let $P\in\mathbb{Z}[x,y]$ be an admissible polynomial.  If the fields $F_j$ associated with the toric points each have $r_j$ real and one pair of complex embeddings, then
\[ \pi m(P) = \sum c_j\,\frac{|\mathrm{disc}(F_j)|^{3/2}}
{\pi^{2r_j+2}} \zeta_{F_j}(2), \qquad c_j\in\mathbb{Q}^*. \]  \label{cortwo} \end{cor}
It is also possible that both corollaries may apply to the same polynomial, as occurs in~(\ref{exthree}).

We now return to the proof of our main result.

\pf  As usual, let $E$ be the field containing the roots of all edge polynomials of~$P$, and choose a parametrization $f(t)=\lambda_1\prod(t-\alpha_r)^{l_r}$ and $g(t)=\lambda_2 \prod(t-\beta_s)^{m_s}$ whose zeros, poles, and coefficients lie in~$E$.  The leading coefficient $\lambda(x)$ of $P(x,y)$ is an edge polynomial of~$P$, hence cyclotomic by assumption.  Therefore $m(\lambda(x))=0$ by Jensen's formula.  Furthermore, the sum
\[ \sum_r l_r \log|\tilde{g}(\alpha_r)| \cdot \mathrm{wind} (\gamma_k,\alpha_r) - \sum_s m_s \log|\tilde{f}(\beta_s)| \cdot \mathrm{wind} (\gamma_k,\beta_s) \]
appearing in~(\ref{tool}) will also vanish using precisely the same reasoning found in the proof of Proposition~\ref{bloch}.  Consequently only the dilogarithm terms contribute to~$m(P)$.  We may also disregard any loops in the set~$S$, since these produce pairs of dilogarithm terms which cancel, as illustrated previously.  These simplifications allow us to rewrite~(\ref{tool}) as
\begin{equation} 2\pi m(P) = \sum_{j=1}^n {\sum_{r,s}}' \pm l_rm_s 
D\!\left(\frac{u_j-\alpha_r}{\beta_s-\alpha_r}\right), \end{equation}
where $(f(u_j),g(u_j))=(\mu,\nu)$ for some toric point of~$P$, and the sign in front of the dilogarithm depends on whether the toric point is an initial point ($+$) or terminal point ($-$) of the corresponding path in the set~$S$ from Theorem~\ref{eval}.  Note that if $(\mu,\nu)$ is singular then it will be an endpoint of more than one path, giving rise to several values of~$u_j$.  In other words, $n$ is equal to twice the number of paths in~$S$, which may be more than the number of toric points of~$P$ if some of them are singular.

Now let $(f_1,g_1)$, \ldots, $(f_w,g_w)$ be the pairs of rational functions obtained by applying the elements of $\mathrm{Gal}(E/\mathbb{Q})$ to each of the coefficients of $f(t)$ and~$g(t)$.  We label their zeros and poles as $\alpha_{rk}$ and $\beta_{sk}$ in the natural manner, so that if $f_k$ is obtained from $f$ via an automorphism $\sigma\in\mathrm{Gal}(E/\mathbb{Q})$, for example, then $\alpha_{rk}=(\alpha_r)^{\sigma}$.  (The factor of $(t-\alpha_{rk})$ in $f_k$ still appears to the power~$l_r$, of course.)  Because $P(x,y)$ has integral coefficients, each pair also provides a parametrization of~$P$.  Let $u_{jk}$ be the corresponding values mapping to the toric points via the $k^{\mathrm{th}}$ parametrization.  Note that $u_{jk}\in K_j=EF_j$ regardless of the parametrization, by Proposition~\ref{edgepoly}.  Applying our evaluation theorem to each pair $(f_k,g_k)$ and summing yields
\[ 2\pi w m(P) = \sum_{j=1}^n \sum_{k=1}^w {\sum_{r,s}}' \pm l_rm_s 
D\!\left(\frac{u_{jk}-\alpha_{rk}}{\beta_{sk}-\alpha_{rk}}\right), \]
where the sign depends only on~$j$.  We next multiply through by a suitable positive integer $N$ so that
Proposition~\ref{bloch} applies to every toric point.  It implies that for a given $j$ and~$k$, the innermost sum over $r,s$ may be written as $D(\theta_{jk})$ for some $\theta_{jk}\in B(K_j)$.  Our formula becomes
\[ 2\pi Nwm(P) = \sum_{j=1}^n \sum_{k=1}^w \pm D(\theta_{jk}). \]
We now claim that the Galois descent property applies to the element $\pm \sum_k \theta_{jk}\in B(K_j)$.  But this is clear from the way the $\theta_{jk}$ were constructed: each $\sigma\in \mathrm{Gal}(K_j/F_j)$ permutes the numbers $\alpha_{rk}$, $\beta_{sk}$, and $u_{jk}$ in a parallel fashion.  Observe that $\sigma$ fixes $\mu_j, \nu_j\in F_j$, so the toric point to which the $u_{jk}$ are mapped does not change.  Hence $\sigma$ simply  permutes the terms of $\sum_k \theta_{jk}$.  We conclude that $\pm \sum_kD(\theta_{jk})=c_jD(\vartheta_j)$ for some $\vartheta_j\in B(F_j)$.  Dividing through by $2Nw$ completes the proof. \fp

\section{Illustrations of the Main Result}

Most of the effort in creating conjectured equalities such as the three given at the start of this paper goes into finding an admissible polynomial to which the corollaries may be applied.  One then simply computes the values of~$m(P)$, the appropriate dilogarithms, and the zeta-functions to a high degree of accuracy and searches for small linear dependencies among them.  Proving the resulting statements is another matter, of course.  To indicate the sorts of steps involved, we now establish~(\ref{exone}), the first of the three original examples.

\begin{prop} If $P(x,y)=y^2+y(x+1)+x^2+x+1$ and $\omega=e^{2\pi i/3}$, then 
\begin{equation} \pi m(P)=2D(i)-{\textstyle\frac{3}{4}}D(\omega). \end{equation} \end{prop}

\begin{figure} \centerline{\includegraphics{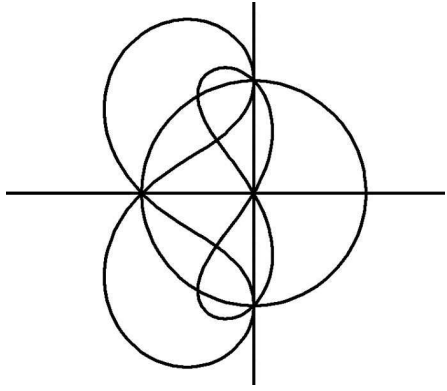}}
\caption{The values of $y$ satisfying $P(e^{i\varphi},y)=0$ for $P(x,y)=y^2+y(x+1)+x^2+x+1$.}
\label{graph3} \end{figure}

\pf The eight toric points of~$P$ are $(-1, i)$, $(i, -1)$, $(i, -i)$, $(\omega,\bar{\omega})$ and their conjugates.  Since $P(x,y)$ is non-singular and tempered, it is admissible.  Therefore only the dilogarithm terms will contribute to~$m(P)$.  A parametrization is provided by Proposition~\ref{param}; we recast it slightly by replacing $t$ by~$-3t$ to obtain
\[ f(t) = \frac{(t+\frac{1}{3})(t-\frac{2}{3})}{t}, \quad g(t) = \omega
\cdot\frac{(t+\frac{1}{3}\bar{\omega})(t-\frac{2}{3}\bar{\omega})}{t}. \]
The set of points $y$ for which $P(x,y)=0$ and $|x|=1$ is shown in Figure~\ref{graph3}.
There are four separate portions which lie outside the unit circle, which give rise to paths $\gamma_1$, \ldots, $\gamma_4$ in the $t$-plane via the above parametrization.  The graph passes through $y=i$
twice, corresponding to $x=-i$ (the outer path) and $x=-1$ (the inner path).  One finds that $f(t)=-1$ and $g(t)=i$ for $t=\frac{1}{3}(-1+\sqrt{3})$; this is the initial point $u_3$ of~$\gamma_3$, the lift of the inner path in our labeling scheme.  In the same manner, we may catalog the remaining initial and terminal points as done in Table~\ref{endpts}.

We label the zeros and poles of $f(t)$ and $g(t)$ by $\alpha_1=-\frac{1}{3}$, $\alpha_2=\frac{2}{3}$, $\alpha_3=0$, $\beta_1=-\frac{1}{3}\bar{\omega}$, $\beta_2=\frac{2}{3}\bar{\omega}$, and $\beta_3=0$. Since $f$ and $g$ have a common pole, there are eight terms in the sum appearing in~(\ref{tool}).  Together with the eight endpoints listed above, we have a total of sixty-four dilogarithms with which to contend, a somewhat daunting prospect.  Fortunately, many of these terms cancel in
pairs, while other values occur several times.  For example, the term corresponding to $\alpha_3$, $\beta_2$, and $u_3$ is
\[ D\left(\frac{\frac{1}{3}(-1+\sqrt{3})-0}{\frac{2}{3}\bar{\omega}-0}
\right) = D\left({\te\frac{1}{4}}(1-\sqrt{3})(1-i\sqrt{3})\right), \]
while the term with $\alpha_3$, $\beta_1$, and $v_4$ is
\[ -D\left(\frac{\frac{1}{3}(-1-\sqrt{3})-0}{-\frac{1}{3}\bar{\omega}-0}
\right) = -D\left({\te\frac{1}{2}}(1+\sqrt{3})(1-i\sqrt{3})\right). \]
But the arguments of these dilogarithms are conjugate reciprocals, and $D(\frac{1}{\bar{z}})=D(z)$, so these two terms will cancel.  In this manner a total of fourteen pairs cancel, while eleven other terms
involve the dilogarithm of a real number, which is zero.  Therefore we are left with twenty-five terms, representing only six distinct values.  Using basic identities such as $D(z)=D(1-\bar{z})=-D\left(\frac{1}{z}\right)$, we may rewrite $2\pi m(P)$ as
\[ 6D(\omega+i\bar{\omega}) + 6D(\bar{\omega}+i\omega) + 6D(i) +
3D(i\omega) + 3D(i\bar{\omega}) - D(1+\omega). \]

\begin{table} \[ \begin{array} {r|c|l} 
\mbox{Endpoints} & \mbox{Path}
& \hspace{2cm} \mbox{Endpoints\ $u_j$,\ $v_j$\ in $t$-plane} \\ \hline 
i,\ -1\ & \gamma_1 & u_1=\frac{1}{6}(1+\sqrt{3}-i(3+\sqrt{3})),
\ v_1=\frac{1}{6}(1+\sqrt{3}+i(3+\sqrt{3})) \rule{0pt}{5mm} \\
-1,\ -i\ & \gamma_2 & u_2=\frac{1}{6}(1-\sqrt{3}-i(3-\sqrt{3})),
\ v_2=\frac{1}{6}(1-\sqrt{3}+i(3-\sqrt{3})) \rule{0pt}{5mm} \\
i,\ \omega\ & \gamma_3 & u_3=\frac{1}{3}(-1+\sqrt{3}),\
v_3=\frac{1}{6}(1-i\sqrt{3})  \rule{0pt}{5mm} \\
\bar{\omega}, -i\ & \gamma_4 & u_4=\frac{1}{3}(-1+i\sqrt{3}),\
v_4=\frac{1}{3}(-1-\sqrt{3}) \rule{0pt}{5mm} 
\end{array} \]
\caption{Coordinates of initial and terminal points.}
\label{endpts} \end{table}

We wish to prove that $2\pi m(P)=4D(i)-{\te\frac{3}{2}} D(\omega)$.  Since $D(1+\omega)= \frac{3}{2}D(\omega)$, it suffices to establish that
\begin{equation} 6D(\omega+i\bar{\omega}) + 6D(\bar{\omega}+i\omega) + 
2D(i) + 3D(i\omega) + 3D(i\bar{\omega}) = 0. \label{suffice} \end{equation}
The time has come for several judicious applications of the five-term identity, which we reproduce below.  It states that
\[ D(a)+D(b)+D\left(\frac{1-b}{a}\right)+ D\left(\frac{a+b-1}{ab}\right)+
 D\left(\frac{1-a}{b}\right)=0 \]
for all $a,b\in\mathbb{C}^*$.  Choosing $a=\omega+i\bar{\omega}$ and 
$b=1-\overline{(\omega+i\bar{\omega})}$ leads to
\[ 2D(\omega+i\bar{\omega}) = D(i) + D(i\omega) + D(1+\omega), \]
while taking $a=\bar{\omega}+i\omega$ and $b=1-\overline{(\bar{\omega}+i\omega)}$ gives
\[ 2D(\bar{\omega}+i\omega) = D(i) + D(i\bar{\omega}) - D(1+\omega). \]
These equations imply that~(\ref{suffice}) is equivalent to $6D(i\omega) + 6D(i\bar{\omega}) + 8D(i) = 0$.  To complete our computation, we employ the Kubert identity for the dilogarithm, discussed in~\cite{Ray}:
\begin{equation} \sum_{k=1}^n D(\xi_n^k z) = \frac{1}{n} D(z^n), \qquad
\xi_n=e^{2\pi i/n}. \label{kubert} \end{equation}
This identity tells us that
\[ 3D(i\omega) + 3D(i\bar{\omega}) + 3D(i) = D(i^3) = -D(i), \]
which confirms the previous equation.  This completes the proof. \fp

We remark that a proof of~(\ref{extwo}) along similar lines exists, in principle.  An attempt to find it came tantalizingly close to succeeding before finally being abandoned.  So we instead apply our main result to the polynomial from~(\ref{exthree}), shown here along with its Newton polygon.
\[  P(x,y) = y^2+y(x^2+x+1)+x^4+x^3+x^2+x+1 \quad\longrightarrow\quad 
\begin{array}{ccccc} 1&&&& \\ 1&0&1&& \\ 1&1&1&1&1 \end{array} \]
This polygon has one interior point, which implies that it generically defines a curve of genus one.  However, the curve has a singularity at $(-1,-1)$ and hence has genus zero.  Its eight toric points are 
$(-\omega,-\bar{\omega})$, $(-\omega,-1)$, $(i,i)$, $(-1,-1)$, and their conjugates.  The fields associated
with the nonsingular toric points are clearly $\mathbb{Q}(\sqrt{-1})$ and~$\mathbb{Q}(\sqrt{-3})$.  To determine the field for the singular point, we expand $P(x-1,y-1)=3x^2-2xy+y^2+\cdots$.   The edge polynomial $3z^2-2z+1$ has roots $\frac{1}{3}(1\pm2i\sqrt{2})$, thus the field is~$\mathbb{Q}(\sqrt{-2})$.
In summary, Corollaries~\ref{corone} and~\ref{cortwo} predict that $\pi m(P)$ will be a rational linear combination of $D(\omega)$, $D(i)$, and a term involving $\zeta_F(2)$ for $F=\mathbb{Q}(\sqrt{-2})$.  We find that to a high degree of accuracy, 
\begin{equation} 5\pi m(P) = 9D(\omega)-\frac{4}{3}D(i)+\frac{16\sqrt{2}}
{\pi^2}\;\zeta_F(2). \end{equation}

As our finale, we examine the polynomial $P(x,y)$ shown below along with its relatively enormous Newton polygon.
\renewcommand{\arraystretch}{.8}
\[ P(x,y) = \frac{x(x+1)^5 - y(y+1)^5}{x-y} \quad\longrightarrow\quad
\begin{array}{cccccc} 1&&&&& \\ 5&1&&&& \\ 10&5&1&&& \\ 10&10&5&1&& \\
5&10&10&5&1& \\ 1&5&10&10&5&1 \end{array} \]
\renewcommand{\arraystretch}{1}
There is a singularity at $(-1,-1)$ of order four, as becomes apparent by expanding
\[ P(x-1,y-1) = -(x^4+x^3y+x^2y^2+xy^3+y^4)+\cdots. \]
The edge polynomial $p(z)=z^4+z^3+z^2+z+1$ defines the field $F=\mathbb{Q}(\xi_5)$ associated with this toric point.  By taking advantage of the form of $P(x,y)$ we are able to find the parametrization
\[ f(t) = -\frac{t^5}{t^5+t^4+t^3+t^2+t+1}, \quad g(t) =
-\frac{1}{t^5+t^4+t^3+t^2+t+1}. \]
The values of $t$ mapping to the singular toric point are $t=\xi_5^{}$, $\bar{\xi}_5^{}$, $\xi_5^2$, and $\bar{\xi}_5^2$.  There are also six other non-singular toric points of the form $(\bar{\xi}_7^k,\xi_7^k,)$ corresponding to~$t=\xi_7^k$, for $k=1$, \ldots, 6.  Each value of $t$ leads to a sum in our expression for $2\pi m(P)$ having the form
\[ \pm 5(D(\omega t)+D(\bar{\omega}t)+D(-\omega t) +D(-\bar{\omega}t)+D(-t)), \]
which simplifies to just $\pm\frac{5}{6}(D(t^6)-6D(t))$ by the Kubert identity.  When $t=\xi_5^k$ this expression further simplifies to $\mp\frac{5}{6}(5D(\xi_5^k))$, while for $t=\xi_7^k$ we obtain $\mp\frac{5}{6}(7D(\xi_7^k))$.  The trickiest part of the whole computation is determining which values of $t$ correspond to initial points of the paths, since a graph is hard to come by.  We find that these values are $t=\xi_5$, $\xi_5^2$, $\bar{\xi}_7$, $\bar{\xi}_7^2$, and~$\bar{\xi}_7^3$.  In summary, we have shown that
\begin{equation} 6\pi m(P) = 35(D(\xi_7^{}) + D(\xi_7^2) +
D(\xi_7^3)) - 25(D(\xi_5^{}) + D(\xi_5^2)). \end{equation}

\end{document}